\documentclass[12pt]{article} 
\usepackage{latexsym} 
\usepackage{amsmath}
\usepackage{amssymb} 
\usepackage{amsfonts}
\usepackage{amscd} 
\usepackage{graphicx}
\usepackage{layout}
\begin{document} 
\newtheorem{Th}{Theorem}[section]
\newtheorem{Cor}{Corollary}[section]
\newtheorem{Prop}{Proposition}[section]
\newtheorem{Lem}{Lemma}[section]
\newtheorem{Def}{Definition}[section]
\newtheorem{Rem}{Remark}[section]
\newtheorem{Ex}{Example}[section]
\newtheorem{stw}{Proposition}[section]


\newcommand{\bet}{\begin{Th}}
\newcommand{\ent}{\stepcounter{Cor}
   \stepcounter{Prop}\stepcounter{Lem}\stepcounter{Def}
   \stepcounter{Rem}\stepcounter{Ex}\end{Th}}


\newcommand{\bec}{\begin{Cor}}
\newcommand{\enc}{\stepcounter{Th}
   \stepcounter{Prop}\stepcounter{Lem}\stepcounter{Def}
   \stepcounter{Rem}\stepcounter{Ex}\end{Cor}}
\newcommand{\bep}{\begin{Prop}}
\newcommand{\enp}{\stepcounter{Th}
   \stepcounter{Cor}\stepcounter{Lem}\stepcounter{Def}
   \stepcounter{Rem}\stepcounter{Ex}\end{Prop}}
\newcommand{\bel}{\begin{Lem}}
\newcommand{\enl}{\stepcounter{Th}
   \stepcounter{Cor}\stepcounter{Prop}\stepcounter{Def}
   \stepcounter{Rem}\stepcounter{Ex}\end{Lem}}
\newcommand{\bef}{\begin{Def}}
\newcommand{\enf}{\stepcounter{Th}
   \stepcounter{Cor}\stepcounter{Prop}\stepcounter{Lem}
   \stepcounter{Rem}\stepcounter{Ex}\end{Def}}
\newcommand{\ber}{\begin{Rem}}
\newcommand{\enr}{
   \stepcounter{Th}\stepcounter{Cor}\stepcounter{Prop}
   \stepcounter{Lem}\stepcounter{Def}\stepcounter{Ex}\end{Rem}}
\newcommand{\bee}{\begin{Ex}}
\newcommand{\ene}{
   \stepcounter{Th}\stepcounter{Cor}\stepcounter{Prop}
   \stepcounter{Lem}\stepcounter{Def}\stepcounter{Rem}\end{Ex}}
\newcommand{\Proof}{\noindent{\it Proof\,}:\ }
\newcommand{\beP}{\Proof}
\newcommand{\enP}{\hfill $\Box$ \par\vspace{5truemm}}

\newcommand{\EE}{\mathbf{E}}
\newcommand{\QQ}{\mathbf{Q}}
\newcommand{\R}{\mathbf{R}}
\newcommand{\C}{\mathbf{C}}
\newcommand{\ZZ}{\mathbf{Z}}
\newcommand{\KK}{\mathbf{K}}
\newcommand{\NN}{\mathbf{N}}
\newcommand{\PP}{\mathbf{P}}
\newcommand{\HH}{\mathbf{H}}
\newcommand{\uuu}{\boldsymbol{u}}
\newcommand{\xxx}{\boldsymbol{x}}
\newcommand{\aaa}{\boldsymbol{a}}
\newcommand{\bbb}{\boldsymbol{b}}
\newcommand{\AAA}{\mathbf{A}}
\newcommand{\BBB}{\mathbf{B}}
\newcommand{\ccc}{\boldsymbol{c}}
\newcommand{\iii}{\boldsymbol{i}}
\newcommand{\jjj}{\boldsymbol{j}}
\newcommand{\kkk}{\boldsymbol{k}}
\newcommand{\rrr}{\boldsymbol{r}}
\newcommand{\FFF}{\boldsymbol{F}}
\newcommand{\yyy}{\boldsymbol{y}}
\newcommand{\ppp}{\boldsymbol{p}}
\newcommand{\qqq}{\boldsymbol{q}}
\newcommand{\nnn}{\boldsymbol{n}}
\newcommand{\vvv}{\boldsymbol{v}}
\newcommand{\eee}{\boldsymbol{e}}
\newcommand{\fff}{\boldsymbol{f}}
\newcommand{\www}{\boldsymbol{w}}
\newcommand{\0}{\boldsymbol{0}}
\newcommand{\lon}{\longrightarrow}
\newcommand{\ga}{\gamma}
\newcommand{\pa}{\partial}
\newcommand{\QED}{\hfill $\Box$}
\newcommand{\id}{{\mbox {\rm id}}}
\newcommand{\Ker}{{\mbox {\rm Ker}}}
\newcommand{\grad}{{\mbox {\rm grad}}}
\newcommand{\ind}{{\mbox {\rm ind}}}
\newcommand{\rot}{{\mbox {\rm rot}}}
\newcommand{\diver}{{\mbox {\rm div}}}
\newcommand{\Gr}{{\mbox {\rm Gr}}}
\newcommand{\LG}{{\mbox {\rm LG}}}
\newcommand{\Diff}{{\mbox {\rm Diff}}}
\newcommand{\Symp}{{\mbox {\rm Symp}}}
\newcommand{\Ct}{{\mbox {\rm Ct}}}
\newcommand{\Uns}{{\mbox {\rm Uns}}}
\newcommand{\rank}{{\mbox {\rm rank}}}
\newcommand{\sign}{{\mbox {\rm sign}}}
\newcommand{\Spin}{{\mbox {\rm Spin}}}
\newcommand{\Sp}{{\mbox {\rm sp}}}
\newcommand{\Int}{{\mbox {\rm Int}}}
\newcommand{\Hom}{{\mbox {\rm Hom}}}
\newcommand{\codim}{{\mbox {\rm codim}}}
\newcommand{\ord}{{\mbox {\rm ord}}}
\newcommand{\Iso}{{\mbox {\rm Iso}}}
\newcommand{\corank}{{\mbox {\rm corank}}}
\def\mod{{\mbox {\rm mod}}}
\newcommand{\pt}{{\mbox {\rm pt}}}
\newcommand{\qed}{\hfill $\Box$ \par}
\newcommand{\spe}{\vspace{0.4truecm}}
\newcommand{\ad}{{\mbox{\rm ad}}}
%
\newenvironment{FRAME}{\begin{trivlist}\item[]
	\hrule
	\hbox to \linewidth\bgroup
		\advance\linewidth by -10pt
		\hsize=\linewidth
		\vrule\hfill
		\vbox\bgroup
			\vskip5pt
			\def\thempfootnote{\arabic{mpfootnote}}
			\begin{minipage}{\linewidth}}{%
			\end{minipage}\vskip5pt
		\egroup\hfill\vrule
	\egroup\hrule
	\end{trivlist}}

\title{
Dependence of vector fields and singular controls
} 

\author{G. Ishikawa
and W. Yukuno}


\date{ }

\maketitle

%
%

\section{Introduction.}

Systems of vector fields on manifolds are important objects 
in various area, say, in geometry, topology, global analysis and control theory. 
Let $X_1, X_2, \dots, X_r$ be a system of 
$C^\infty$ vector fields over a $C^\infty$ manifold $M$ of dimension $n$, with $r \leq n$. 
We consider the control system 
$$
\dot{x} = \sum_{i=1}^r u_i X_i(x), 
$$
which is called a driftless control-affine system, generated by $X = (X_1, X_2, \dots, X_r)$ 
with the control parameters $u = (u_1, \dots, u_r)$ in an open set $\Omega \subset \R^r$. 

In general, even generically, the {\it dependence locus} 
$$
\Sigma := \{ x \in M \mid X_1(x), X_2(x), \dots, X_r(x) {\mbox{\rm \ are linearly dependent}} \}
$$
is need not be void. 
We show, in this paper, an example providing a significance in geometric control theory 
of the existence of the dependence locus of the system. 
We show, in particular, the generic appearance of non-trivial 
singular trajectories embedded in the dependence locus. 

An absolutely continuous curve $x : I \to M$ defined on a closed interval $I$ is called an 
{\it $X$-trajectory} if it is a solution curve of the differential equation 
$\dot{x} = \sum_{i=1}^r u_i(t) X_i(x)$ for 
some essentially bounded measurable function $u : I \to \Omega$ which is called a control. 
Note that the control is not unique where $X_1, X_2, \dots, X_r$ are linearly dependent. 
A control $u : I \to \Omega$ and an initial point determine the trajectory $x = x(t)$ if it exists. 
Suppose a control $u : I \to \Omega$ defines an $X$-trajectory $x : I \to M$ on $I$ with an initial point 
$x_0 = x(a), I = [a, b]$. Then any $\widetilde{u} : I \to \Omega$ in a neighborhood ${\mathcal C}$ of $u$,  
in the $L^\infty$ topology, defines uniquely an $X$-trajectory $\widetilde{x} : I \to M$ on the same interval $I$ 
with the initial point $x_0 = \widetilde{x}(a)$, so it defines the endpoint 
$\widetilde{x}(b)$. Thus we have the endpoint mapping ${\mathcal C} \to M$, $\widetilde{u} 
\mapsto \widetilde{b}$. Then $u$ is called a {\it singular control} if it is a singular point of the endpoint mapping, 
and the corresponding trajectory is called a {\it singular $X$-trajectory} (see \cite{AS}\cite{BC}\cite{CJT}\cite{CJT2}). 

An $X$-trajectory $x : I \to M$ 
is called {\it non-trivial} (resp. {\it dependent}) if
for any subinterval $J \subset I, x\vert_J : J \to M$ is not constant (resp. 
$x(I) \subset \Sigma$). 

\

\noindent
{\bf Theorem.} 
{\it 
Let $M$ be a $C^\infty$ manifold of dimension $3$. Then there exists 
a non-empty open subset ${\mathcal U}$ in the set of systems of triple $C^\infty$ vector fields 
$X = (X_1, X_2, X_3)$ over $M$ with the $C^\infty$ topology 
such that, for any $X \in {\mathcal U}$, there exists 
a one parameter family of non-trivial singular $X$-trajectories 
which are embedded in the dependent locus $\Sigma$ of $X$ and 
lifted to essentially unique singular controls. 
}

\

The above Theorem means that there exists an example of a system with 
the property stated in Theorem such that any small perturbation of the system 
for $C^\infty$ topology enjoys the same property. Therefore, in general, 
it is impossible to exclude generically the appearance of non-trivial singular trajectories 
embedded in the dependent locus. 


\


The generic properties of driftless control-affine systems 
are established in \cite{CJT} 
when the dependent locus is void. Moreover in \cite{CJT2} the theory was extended 
to the general case where the dependent locus may not be void. 
However a theorem (Theorem 2.13 of \cite{CJT2}) used there, which says that
any $X$ trajectory has essentially zero velocity on dependent locus 
for a generic $X$ with $r \leq n$, is incorrect and requires more conditions. 
We provide an example to the approach to the general theory based 
on the existence of dependent locus. 

Regarding the existence of non-trivial dependent singular trajectories, 
the second author reformulates a theorem of \cite{CJT2} (Theorem 2.17 of \cite{CJT2}) 
and gives its proof in \cite{Yukuno}. 

\

We concern with only small perturbations for $C^\infty$ topology. 
However if $M$ is not parallelizable, i.e. $TM$ is not trivial, 
for instance if $M = S^3$, any system $X = (X_1, X_2, X_3)$ on $M$ has non-void 
dependence locus $\Sigma$. 
Thus the phenomena of our Theorem 
are unavoidable by any means. 

On the removability of dependence locus by large perturbations, 
see, for instance, the famous topological studies \cite{AD}\cite{Thomas}. 

%

\section{Proof of Theorem.}

Let $M$ be a manifold of dimension $3$. 
We consider a generic distribution $D = \langle X_1, X_2, X_3 \rangle \subset TM$ with 
singularities. In fact we impose on $X = (X_1, X_2, X_3)$ the following generic conditions: 

(1) The dependence locus $\Sigma$ is a smooth surface. Moreover, 

(2) Outside of a curve $\gamma \subset \Sigma$, 
$D_x$ is transverse to $T_x\Sigma$ ($x \in \Sigma \setminus \gamma$). 

The above conditions are achieved by a transversality of $1$-jets of $X$. 
Moreover for any $M$, there exists a system $X$ satisfying the above transversality condition 
with $\Sigma \not= \emptyset$. 

Take any $x_0 \in \Sigma \setminus \gamma$. Then there exists a system of local coordinates 
$x_1, x_2, x_3$ around $x_0$ on an open set $U \subset M$ such that $\Sigma = \{ x_1 = 0\}$ 
and $D$ is generated by 
$$
\left\{
\begin{array}{cccccc}
X_1 & = & \frac{\pa}{\pa x_1} & & + & P\ \frac{\pa}{\pa x_3}
\smallskip
\\
X_2 & = & & \frac{\pa}{\pa x_2} & + & Q\ \frac{\pa}{\pa x_3}
\smallskip
\\
X_3 & = & & & & x_1 \ \frac{\pa}{\pa x_3}
\end{array}
\right., 
$$
where $P = P(x_2, x_3), Q = Q(x_2, x_3)$. Then we have 
\begin{equation*}
\begin{split}
D = & \ \{ (x_1, x_2, x_3; u_1, u_2,  P(x_2, x_3)u_1 + Q(x_2, x_3)u_2 + x_1u_3) 
\\
 & 
\hspace{7truecm}
\mid (x_1, x_2, x_3) \in U, u_1, u_2, u_3 \in \R\}. 
\end{split}
\end{equation*}
and the dependence locus $\Sigma = \{ x_1 = 0\}$. Then we have 
$$
D \cap T\Sigma = \{ (0, x_2, x_3; 0, u_2, Q(x_2, x_3)u_2) \}, 
$$
which is a line field on $\Sigma$ defined by $dx_3 - Q(x_2, x_3)dx_2 = 0$. 

The Hamiltonian function $H : T^*\R^3 \times\R^3 \to \R$ of the system is given by 
$$
\begin{array}{rcl}
H & = & \langle p, u_1X_1 + u_2X_2 + u_3X_3\rangle 
\\
 & = & u_1(p_1 + Pp_3) + u_2(p_2 + Qp_3) + u_3x_1p_3 
\\
 & = & u_1p_1 + u_2p_2 + (u_1P + u_2Q + u_3x_1)p_3
\end{array}
$$
The constrained Hamiltonian system is given by 
$$
\left\{
\begin{array}{l}
\dot{x}_1 = u_1, \ \dot{x}_2 = u_2, \ \dot{x}_3 = u_1P + u_2Q + u_3x_1, 
\\
\dot{p}_1 = - u_3p_3, \ \dot{p}_2 = - (u_1P_{x_2} + u_2Q_{x_2})p_3, \ 
\dot{p}_3 = - (u_1P_{x_3} + u_2Q_{x_3})p_3, 
\\
p_1 + Pp_3 = 0, \ p_2 + Qp_3 = 0, \ x_1p_3 = 0, \quad p = (p_1, p_2, p_3) \not= 0. 
\end{array}
\right.
$$
If $x_1 \not= 0$, 
 then $p = 0$. 
 Therefore we suppose $x_1 = 0$. Then $u_1 = 0$. Then the system is reduced to the system on $\Sigma$: 
$$
\left\{
\begin{array}{l}
\dot{x}_2 = u_2, \ \dot{x}_3 = u_2Q, 
\\
 \dot{p}_2 = - u_2Q_{x_2}p_3, \  \dot{p}_3 = - u_2Q_{x_3}p_3, 
\\
p_2 = - Qp_3, \ p_3 \not= 0, 
\end{array}
\right.
$$
with additional conditions $p_1 = - Pp_3, \ u_3 = - \dot{p}_1/p_3, \ x_1 = 0$. 

Take any solution $x_2 = x_2(t), x_3 = x_3(t)$ of the equation $dx_3 - Q(x_2, x_3)dx_2 = 0$. 
Then we have 
$$
p_3(t)  = a \exp \int \left(- \dot{x}_2(t)Q_{x_3}(x_2(t), x_3(t))\right) dt, \ a \not= 0, 
$$
and therefore $\dot{p}_2(t)$ and $p_1(t) = - P(x_2(t), x_3(t))p_3(t)$ are determined. 
Then we have 
$$
(p_2 + Qp_3)' = \dot{p}_2 + Q_{x_2}\dot{x}_2p_3 + Q_{x_3}\dot{x}_3p_3 + Q\dot{p}_3 
= 0. 
$$
If we choose any initial value of $p_2$, then we have $p_2 = - Qp_3$. 
Therefore, for any solution curve $(x_2(t), x_3(t))$ 
of the equation $dx_3 - Q(x_2(t), x_3(t))dx_2 = 0$, 
the curve $(0, x_2(t), x_3(t))$ on the dependence locus $\Sigma$ is a singular trajectory 
of the control system. Moreover the corresponding singular control is unique to it. 
\QED

%
%
%
%

\

\begin{flushleft}
Goo ISHIKAWA, \\
Department of Mathematics, Hokkaido University, 
Sapporo 060-0810, Japan. \\
e-mail : ishikawa@math.sci.hokudai.ac.jp \\

\

Wataru YUKUNO, \\
Department of Mathematics, Hokkaido University, 
Sapporo 060-0810, Japan. \\
e-mail : yukuwata@math.sci.hokudai.ac.jp

\end{flushleft}

\end{document}